\documentclass[12pt,reqno]{article}
mm \textheight=8.9in

\usepackage{amssymb}
\usepackage{amsmath}
\usepackage{amsthm}
\usepackage{amssymb}
\usepackage{amsmath}
\usepackage{amsthm}
\usepackage{graphics}
\usepackage{color}
\newcommand{\tw}[3]{{$#1$}${\,\scriptscriptstyle {#2}}\atop\raise9pt\hbox{$\scriptstyle\tp$} ${$#3$}}
\newcommand{\st}[1]{\mbox{${\,\scriptscriptstyle {#1}}\atop\raise5.5pt\hbox{$*$}$}}
\newcommand{\btr}{\raise1.2pt\hbox{$\scriptstyle\blacktriangleright$}\hspace{2pt}}
\newcommand{\braid}[3]{{#1}$\lower4pt\hbox{$\tp\atop\raise4pt
            \hbox{$\scriptscriptstyle {#2} $}$}${#3}}

\newcommand{\id}{\mathrm{id}}
\newcommand{\im}{\mathrm{im}\:}

\newcommand{\A}{\mathcal{A}}

\newcommand{\Dc}{\mathcal{D}}

\newcommand{\Yg}{\mathrm{Y}}
\newcommand{\Ug}{\mathrm{U}}
\newcommand{\Xg}{\mathrm{X}}

\newcommand{\Sg}{\mathfrak{S}}

\newcommand{\Mg}{\mathfrak{M}}

\newcommand{\Ha}{\mathcal{H}}
\newcommand{\Ru}{\mathcal{R}}

\newcommand{\Sc}{\mathcal{S}}

\newcommand{\C}{\mathbb{C}}
\newcommand{\Z}{\mathbb{Z}}

\newcommand{\tp}{\otimes}
\newcommand{\vt}{\vartheta}
\newcommand{\zt}{\zeta}
\newcommand{\V}{V}

\newcommand{\dt}{\delta}

\newcommand{\la}{\lambda}

\newcommand{\End}{\mathrm{End}}

\newcommand{\Tr}{\mathrm{Tr}}

\newcommand{\btl}{\mbox{\raise1.1pt\hbox{$\scriptstyle\blacktriangleleft$}}}

\newcommand{\g}{\mathfrak{g}}

\newcommand{\s}{\mathfrak{s}}
\renewcommand{\o}{\mathfrak{o}}

\newcommand{\nn}{\nonumber}

\newcommand{\p}{\mathfrak{p}}
\renewcommand{\l}{\mathfrak{l}}

\newcommand{\si}{\sigma}
\newcommand{\al}{\alpha}

\newcommand{\bt}{\beta}

\newcommand{\be}{\begin{eqnarray}}
\newcommand{\ee}{\end{eqnarray}}

\newtheorem{thm}{Theorem}
\newtheorem{propn}[thm]{Proposition}
\newtheorem{lemma}[thm]{Lemma}

\theoremstyle{definition}

\begin{document}
\title{Irreducibility of fusion modules over twisted Yangians at generic point}
\author{A. I. Mudrov}
\date{}
\maketitle
\begin{center}
{
Department of Mathematics, University of York, YO10 5DD, UK,\\
e-mail: aim501@york.ac.uk}\\
 St.-Petersburg Department of Steklov Mathematical Institute,
\\
Fontanka 27, 191011 St.-Petersburg, Russia,\\
 \texttt{mudrov@pdmi.ras.ru}
\end{center}
\begin{abstract}
With any skew Young diagram one can associate a one parameter family
of "elementary" modules over the Yangian $\Yg(\g\l_N)$.
Consider the twisted Yangian $\Yg(\g_N)\subset \Yg(\g\l_N)$ associated
with a classical matrix Lie algebra $\g_N\subset\g\l_N$.
Regard the tensor product of elementary Yangian modules as
a module  over $\Yg(\g_N)$ by restriction. We prove its irreducibility
for generic values of the parameters.
\end{abstract}
{\small \underline{Key words}:
reflection equation, twisted Yangian, universal S-matrix, fusion modules.\\
\underline{AMS classification codes}: 81R50,81R10, 20G42.}
\section{Introduction}
We study irreducibility of a special class of modules
over twisted Yangians. Twisted Yangians are associated with a certain type of reflection equation
and appear in integrable spin chains with boundaries,
alongside with the Yangian $\Yg(\g\l_N)$, which is an RTT-type algebra.
The twisted Yangian is a subalgebra
in  $\Yg(\g\l_N)$, and the modules in question
are obtained from the so called fusion  $\Yg(\g\l_N)$-modules by restriction.
Our goal is to show that generic fusion module over the twisted Yangian is irreducible,
thus extending a result of \cite{NT} on fusion $\Yg(\g\l_N)$-modules.

More specifically, consider the Lie algebra $\g_N\subset \g\l_N$ of matrices preserving a non-degenerate symmetric
or skew symmetric inner product in $\C^N$. One can relate to
$\g_N$ a a certain subalgebra  $\Yg(\g_N)$ in $\Yg(\g\l_N)$ called the twisted  or
Olshanski Yangian.
The defining representation of $\g\l_N$ in $\C^N$ induces
a one parameter family $V(z)$ of representations of  $\Yg(\g\l_N)$ through
the composition $\Yg(\g\l_N)\to \Yg(\g\l_N)\to \Ug(\g\l_N)$, where the right arrow
is the canonical evaluation homomorphism while the left one is the so called shift automorphism
depending on $z$.
Fusion $\Yg(\g\l_N)$-modules are obtained from $V(z)$ via tensor product (using the Hopf algebra structure of the Yangian) and taking subquotients.
Among them there are shifted evaluation modules, that are factored
through the evaluation homomorphism $\Yg(\g\l_N)\to \Ug(\g\l_N)$
(this is not the case in general, because the evaluation homomorphism is not a coalgebra map).
Such modules are parameterized with Young diagrams.
A more general class of modules involves
the so called skew Young diagrams (pairs of subordinate Young diagrams)
and called elementary in \cite{NT}. Every elementary module is included into a one parameter family with
the same underlying vector space. In the present paper
we understand by fusion $\Yg(\g\l_N)$-module
the tensor product of elementary ones. Thus, a fusion $\Yg(\g\l_N)$-module is included into a family
parameterized by a point of an affine space.

As $\Yg(\g_N)$ is a subalgebra in $\Yg(\g\l_N)$, every representation of the Yangian becomes
a representation of the twisted Yangian.
It has been shown in \cite{NT}
that generic fusion $\Yg(\g\l_N)$-module is irreducible. We prove that
the fusion modules  are
as well irreducible over the twisted Yangians, except for an algebraic set in
the  parameter space.

To prove irreducibility, we employ a generalization of the method used
in \cite{NT} (see also \cite{KMT} and \cite{MT}) for   $\Yg(\g\l_N)$-modules
and  in \cite{MN} for special $\Yg(\g_N)$-modules.
The method is based on the following  observation. Suppose that $\A$ is an associative
algebra and $(Z,\rho)$ is a finite dimensional  representation of $\A$. Let $W$ be an auxiliary
vector space. Suppose there exists a linear operator $\Phi_0\colon \End(W)\to \End(Z)$
whose image lies in $\rho(\A)$. If $\Phi_0$ is onto, then $Z$ is irreducible.
The operator $\Phi_0$ is taken in \cite{NT} to be the leading
Laurent term of the Yangian R-matrix $R_{W,Z}$ at a certain point.
Here $R_{W,Z}$ is considered as a rational function of a complex variable parameterizing the module $W$.
In the present paper we use the so called universal S-matrix, to obtain $\Phi_0$.

The universal S-matrix is an element  $\Sc\in \Ha\tp \Yg(\g_N)$, where $\Ha$ is
a certain pseudo quasitriangular bialgebra equipped with an involution $\tau$, see \cite{M}. The prefix "pseudo" means that the
universal R-matrix $\Ru$ of $\Ha$ is not invertible. The involution $\tau$ is an anti-algebra and a coalgebra
automorphism. The matrix $\Sc$ obeys the equality
$(\Delta\tp \id)(\Sc)=\Sc_1\Ru'\Sc_2$ with $\Ru':=(\tau\tp \id)(\Ru)$. It can  be
constructed out of the generators matrix  of the algebra $\Yg(\g_N)$ through
a fusion procedure facilitated by the above equality.

The algebra $\Ha$ admits a non-trivial homomorphism
from the Yangian $\Yg(\g\l_N)$ and hence from the twisted Yangian $\Yg(\g_N)$.
Under that homomorphism, every fusion module is restricted from a representation of $\Ha$.

Given a fusion  module $Z=Z(z_1,\ldots, z_\ell)$ we take  the $\Ha$-module $Z(z_1+\zt,\ldots, z_\ell+\zt)$ for the auxiliary space $W$.
Let $\rho$ denote the representation homomorphism $\Yg(\g_N)\to \End(Z)$ and let  $\Sc_{W,Z}$ be the image of $\Sc$ in $\End(W)\tp \End(Z)$.
The right tensor component of the matrix $\Sc_{W,Z}$ lies in $\rho\bigl(\Yg(\g_N)\bigr)\subset \End(Z)$.
Expanding  $\Sc_{W,Z}$ in the Laurent series around $\zt=0$ we take the leading term for the operator $\Phi_0$.
We prove that $\Phi_0$ implements a surjective mapping $\End(Z)\to \End(Z)$  for almost all values of the
parameters $(z_1,\ldots, z_\ell)$.
In particular, the module $Z(z_1,\ldots, z_\ell)$ is irreducible if all the parameters $z_i$ are not half-integers
and $z_i\pm z_j\not \in \Z$ for all pairwise distinct $i,j$.

The paper is organized as follows.
Section \ref{YtY} contains definitions and basic properties of the Yangian $\Yg(\g\l_N)$ and the twisted Yangians $\Yg(\g_N)$.
In Section \ref{USM} we introduce, following \cite{M}, the auxiliary bialgebra $\Ha$ and
the universal S-matrix $\Sc\in \Ha\tp \Yg(\g_N)$. Section \ref{Sec:H-mod} contains
general facts about $\Ha$-modules. In Section \ref{sYd} we consider modules associated with
skew diagrams.  Although the constructions of this section are standard for $\Yg(\g\l_N)$-modules,
their extension to $\Ha$-modules is important for our approach.
In Section \ref{Cd} we study the involutive operation on skew diagrams obtained by
$180^\circ$-rotation around the diagram center. This operation gives rise
to an operation on $\Ha$-modules, which is studied in Section \ref{CHm}.
Section \ref{Sec:factorization} contains an auxiliary information on factorization
of R- and S-matrices.
The main result of the paper is proved in Section \ref{Irreducibility_gen}.

\vspace{0.2in}
\noindent
{\bf \large Acknowledgements}

The author is grateful to Maxim Nazarov for formulating the problem and for valuable discussions.
This research is supported by the EPSRC grant C511166
and partly supported by the RFBR grant 06-01-00451.

\section{Yangian  and twisted Yangians}
\label{YtY}
In this section we recall the definition of the Yangian $\Yg(\g\l_N)$ and the twisted Yangian $\Yg(\g_N)$,
following \cite{MNO}.

Fix a positive integer $N>1$. In the associative algebra $\End(\C^N)$ of linear operators
on the vector space $\C^N$ choose the standard matrix basis $\{e_{ij}\}_{i,j=1}^N \subset \End(\C^N)$
with the multiplication  $e_{ij}e_{mn}=\dt_{jm}e_{in}$. Denote by $\g\l_N$ the general linear Lie algebra $\End(\C^N)$.

Fix an involutive anti-algebra isomorphism $A\mapsto A^t$ of $\End(\C^N)$. In is factored
through the matrix transposition and the subsequent similarity transformation $A\mapsto gAg^{-1}$ with
a non-degenerate symmetric or skew-symmetric matrix $g\in \End(\C^N)$.
Denote by $\g_N\subset \g\l_N$ the Lie algebra of matrices obeying the condition $A^t=-A$.
If $g$ is symmetric, then $\g_N=\s\o_N$; otherwise  $\g_N=\s\p_N$.

Let $P$ denote the flip of the tensor factors in $\C^N\tp \C^N$; in the standard basis,
$P=\sum_{i,j=1}^N e_{ij}\tp e_{ji}$.
The element  $R(u,v):=u-v-P\in \End^{\tp 2}(\C^N)$ satisfies the
Yang-Baxter equaiton
$$
R_{12}(u_1,u_2)R_{13}(u_1,u_3)R_{23}(u_2,u_3)=
R_{23}(u_2,u_3)R_{13}(u_1,u_3)R_{12}(u_1,u_2)
$$
and is called the Yang R-matrix.
Here the subscripts indicate the embedding of the algebra $\End^{\tp 2}(\C^N)$ in $\End^{\tp 3}(\C^N)$ in the usual way.

The Yangian $\Yg(\g\l_N)$ is an associative algebra with unit, a deformation of the universal enveloping of the current Lie algebra $\g\l_N[u]$.
It is generated by the infinite set of elements $T^{(k)}_{ij}$, where $i,j=1,\ldots,N$ and $k=1,2\ldots$.
To write down the commutation relations on $T^{(k)}_{ij}$, it is convenient to organize them in the matrix
$T(u)=\sum_{i,j=1}^N e_{ij}\tp T_{ij}(u)$, where $T_{ij}(u):=\sum_{k=0}^{\infty}T^{(k)}_{ij}u^{-k}\in \Yg(\g\l_N)[[u^{-1}]]$.
Here $T^{(0)}_{ij}:=\delta_{ij}$.
Then the defining relations of the algebra  $\Yg(\g\l_N)$ read
$$
R_{12}(u,v)T_{1}(u)T_{2}(v)=T_{2}(v)T_{1}(u)R_{12}(u,v).
$$

The extended twisted Yangian $\Xg(\g_N)$ is an associative unital algebra generated by the infinite set of elements $S^{(k)}_{ij}$, where $i,j=1,\ldots,N$ and $k=1,2\ldots$.
Similarly to the Yangian,  organize the generators $S^{(k)}_{ij}$ in the matrix
$S(u)=\sum_{i,j=1}^N e_{ij}\tp S_{ij}(u)$, where $S_{ij}(u):=\sum_{k=0}^{\infty}S^{(k)}_{ij}u^{-k}\in \Xg(\g_N)[[u^{-1}]]$ and
$S^{(0)}_{ij}:=\delta_{ij}$.
Then the defining relations of the algebra  $\Xg(\g_N)$ read
$$
R_{12}(u,v)S_{1}(u)R'_{12}(u,v)S_{2}(v)=S_{2}(v)R'_{12}(u,v)S_{1}(u)R_{12}(u,v).
$$
Here $R'_{12}(u,v):=-(u+v+Q)$ with $Q:=\sum_{i,j=1}^{N} e^t_{ij}\tp e_{ji}\in \End^{\tp 2}(\C^N)$.

The extended twisted Yangian is a right comodule algebra over the Yangian under the coaction
\be
\label{coaction}
S_{ij}(u)\mapsto \sum_{k,l=1}^NS_{kl}(u)\tp T^t_{ik}(-u)T_{lj}(u),
\ee
where the tensor product is taken over $\C[[u^{-1}]]$.
The assignment
$S(u)\mapsto T^t(-u)T(u)$
defines an algebra homomorphism $\Xg(\g_N)\mapsto \Yg(\g\l_N)$.
The image of this homomorphism is called the twisted Yangian
and denoted by $\Yg(\g_N)$. It is a deformation of the universal enveloping of the twisted current Lie algebra $\g\l_N^t[u]$.
\section{The universal S-matrix}
\label{USM}
The universal S-matrix of the twisted Yangian is an element\footnote{Completion of the tensor product based on a natural filtration in $\Ha$ is understood.}
$\Sc\in \Ha\tp \Xg(\g_N)$, where $\Ha$ is a pseudo quasitriangular
bialgebra with involution to be described later on.
The matrix $\Sc$ satisfies the characteristic equation
\be
\label{char_eq}
(\Delta\tp \id)(\Sc)=\Sc_1\Ru_{12}'\Sc_2,
\ee
where $\Ru':=(\tau\tp \id)(\Ru)$ is obtained from the universal R-matrix $\Ru\in \Ha\tp \Ha$ via an involutive
anti-algebra and coalgebra automorphism $\tau\colon\Ha\to \Ha$. Under the assumption
$(\Ru')_{21}=\Ru'$, the universal S-matrix satisfies the reflection equation
$$
\Ru_{12}\Sc_1\Ru_{12}'\Sc_2=\Sc_2\Ru_{12}'\Sc_1\Ru_{12}
$$
in $\Ha\tp \Ha\tp \Xg(\g_N)$ and plays an important role in our study.

To describe the bialgebra $\Ha$, denote by $\Mg$ the ring of polynomials in the variables $u$ and $u^{-1}$ with coefficients in $\End(\C^N)$
(the loop ring).
Consider the infinite direct sum of algebras $\Sigma\Mg:=\oplus_{k=0}^\infty \Mg^{\tp k}$, where the zeroth
summand is the ground field $\C$. The elements of $\Sigma\Mg$ are infinite formal sums
$\sum_{k=0}^\infty h^{(k)}$, where $h^{(k)}\in\Mg^{\tp k}$.
Define the comultiplication  on each $\Mg^{\tp k}\subset \Sigma\Mg$
as $\sum_{i+j=k}\pi_{ij}$, where $\pi_{ij}\colon \Mg^{\tp (i+j)}\to \Mg^{\tp i}\tp \Mg^{\tp j}$
is the tautological isomorphism. The comultiplication extends over $\Sigma\Mg$ and makes it
a bialgebra, with the counit being the projection $h\mapsto h^{(0)}\in \C$.

Define $\Ha\subset \Sigma\Mg$ to be the subalgebra
$$
\Ha:=\{h\in \Sigma\Mg\>|\> R_{i,i+1}h^{(k)}=\sigma_{i,i+1}(h^{(k)})R_{i,i+1}, \quad i=1,\ldots k-1\}.
$$
Here $\si_{i,i+1}$ is the flip of the $i$ and $i+1$ tensor factors in $\Mg^{\tp k}$.
One can check that $\Ha$ is a sub-bialgebra in  $\Sigma\Mg$. Moreover, it is pseudo quasitriangular
with the (non-invertible) universal R-matrix
\be
\label{R-fusion}
\Ru:=\sum_{k,m=0}^\infty \Ru^{(k),(m)},
\quad
\Ru^{(k),(m)}:=
\left\{
\begin{array}{ccc}
1^{(k)}\tp 1^{(m)}, &km=0,\\
\prod_{i=1}^n\prod_{j=n}^1 R_{i,j}
, &km\not =0.
\end{array}
\right.
\ee
Here $1^{(k)}$ is the unit in $\Mg^{\tp k}$, so that the sum $\sum_{k=0}^\infty 1^{(k)}$ is the unit in $\Sigma\Mg$.
The indices $i,j$ in the ordered product defining  $\Ru^{(k),(m)}$ refer, respectively,  to the first
and the second tensor factors of $\Mg^{\tp k}\tp \Mg^{\tp m}$. Here and further on we adopt the
following convention for  ordered products. The left-most and, respectively, the right-most factors correspond
to the initial  and final values  of the product index.
Alternatively, when taking product over ordered sets we indicate the direction of increment with arrows.

Denote by $\tau$ the involutive anti-automorphism of the ring $\Mg$ defined by
$A\to A^t$ and $u\mapsto -u$, for all $A\in \End(\C^N)$. It naturally extends to any tensor power $\Mg^{\tp k}$ and
hence over the entire $\Sigma\Mg$. It is easy to see that  $\tau$ preserves the comultiplication in $\Sigma\Mg$ and
restricts to the sub-bialgebra $\Ha\subset \Sigma\Mg$. Thus $\Ha$ is equipped with an involutive anti-algebra and
coalgebra automorphism $\tau$. The element $\Ru'=(\tau\tp \id)(\Ru)$ can be constructed out of $R'(u,v)$ via a fusion procedure,
or directly applying $\tau$ to $\Ru$ in the formula (\ref{R-fusion}).
Note that $\Ru'=(\Ru')_{21}$.

It is also  convenient to introduce the matrices
\be
\breve{R}(u,v)&:=&1-\frac{P}{u-v}=\frac{R(u,v)}{u-v},\nn\\
\breve{R}'(u,v)&:=&1+\frac{Q}{u+v}=-\frac{R'(u,v)}{u+v}.\nn
\ee
These matrices can be  viewed as elements of the corresponding "localization" of the ring $\Mg\tp \Mg$.
Using the same fusion procedure as for $\Ru$ and $\breve{\Ru}$, one can define the universal elements $\breve{\Ru}$ and $\breve{\Ru}'=(\tau\tp \id)(\breve{\Ru})$
belonging to a ceratin extension of $\Ha\tp \Ha$. Note that $\Ru$ and $\Ru'$ differ from $\breve{\Ru}$ and, respectively, $\breve{\Ru}'$
by certain central multipliers.

\section{On representations of $\Ha$}
\label{Sec:H-mod}
The Yangian $\Yg(\g\l_N)$ admits a non-trivial homomorphism (defining representation) to the algebra $\Mg$ via the
the assignment
$
T^{k}_{ij}\mapsto - v^{k-1}e_{ji}.
$
In the matrix form this  assignment reads
\be
\label{defining_hom}
T(u)\mapsto 1-\sum_{k=0}^\infty\frac{v^k}{u^{k+1}}P=\breve{R}(u,v).
\ee
The defining representation (\ref{defining_hom}) extends to a bialgebra homomorphism  $\Yg(\g\l_N)\to \Ha$, which
induces a homomorphism $\Xg(\g_N)\to \Ha$ of right $\Ha$-comodule algebras. Hence every $\Ha$-module becomes
automatically a $\Yg(\g\l_N)$- and therefore an $\Xg(\g_N)$-module.
By construction, every such a representation of $\Xg(\g_N)$ factors through a representation of $\Yg(\g_N)$.

Let $\chi\colon \Xg(\g_N)\to \C$ denote the character defined by the assignment $S_{ij}(u)\mapsto \delta_{ij}$.
The homomorphism $\Xg(\g_N)\to \Yg(\g\l_N)$ is factorized
in the composition
$$
\Xg(\g_N)\to \Xg(\g_N)\tp \Yg(\g\l_N)\stackrel{\chi\tp \id}{\longrightarrow} \Yg(\g\l_N),
$$
where the left arrow is the coaction (\ref{coaction}).
The image of the universal S-matrix in $\Ha\tp \Ha$ is equal to
$\breve{\Ru}'_{12}\Sc^\chi_1\breve{\Ru}_{12}$, where $\Sc^\chi:=(\id\tp \chi)(\Sc) \in \Ha$.
The element $\breve {\Ru}\in \Ha\tp \Ha$ has been introduced  in  the previous section.

Let $\pi^{(\ell)}$ denote the projection $\Ha\to \Ha\cap\Mg^{\tp \ell}$.
The  projection $\pi:=\pi^{(1)}$ can be considered as a representation of $\Ha$
in the vector space $V(v):=\C^N[v,v^{-1}]$. Consider the tensor product
$V(v_1)\tp \ldots \tp V(v_\ell)$ as a module over the bialgebra $\Ha$ constructed
through the $\ell$-fold comultiplication $\Delta^{(\ell)}\colon \Ha\to \Ha^{\tp \ell}$.
One can easily see that $\pi^{\tp \ell}\circ \Delta^{(\ell)}=\pi^{(\ell)}$. In other words,
the representation is non-zero only on the subalgebra $\Ha\cap\Mg^{\tp \ell}$, where it is tautological.
It vanishes on every subspace $\Ha\cap\Mg^{\tp k}$ with $k\not=\ell$.

In fact, we shall consider certain extensions of the vector space $V(v_1)\tp \ldots \tp V(v_\ell)$ corresponding
to extensions of the algebra $\C[v_i,v^{-1}_i]$. We assume that wherever necessary
without noticing.

Consider the element $\Ru^{(\ell)}:=\overrightarrow{\prod}_{p<q}\Ru_{p,q}$ in $\Ha^{\tp \ell}$ intertwining the $\ell$-fold comutliplication
with its opposite,
$\Ru^{(\ell)}\Delta^{(\ell)}(h)=\Delta^{(\ell)}_{op}(h)\Ru^{(\ell)}$ for all $h\in \Ha$.
Here the product is taken over lexicographically ordered pairs $1\leqslant p<q\leqslant \ell$; the direction of increment is always
from left to right.

Introduce the operator
$
R^{(\ell)}=\pi^{\tp \ell}(\Ru^{(\ell)})\in \Mg^{\tp\ell}
$
intertwining the representation $V(v_1)\tp \ldots \tp V(v_\ell)$
with the representation $V(v_1)\tilde \tp \ldots \tilde \tp V(v_\ell)$ defined through the opposite comultiplication $\Delta_{op}$.
Hence the image of $R^{(\ell)}$ is a $\Ha$-submodule in $V(v_1)\tilde \tp \ldots \tilde \tp V(v_\ell)$.
Note that for every $h\in  \Ha\cap\Mg^{\tp \ell}$ one has
\be
\label{intertwining}
R^{(\ell)} h=\si_\ell(h)R^{(\ell)},
\ee
where $\si_\ell$ is the involution reversing the order of the tensor factors in $\Mg^{\tp\ell}$.

For any pair of $\Ha$-modules $U$ and $W$  we denote by $\Ru_{U,W}$ and  $\Ru'_{U,W}$
the images of, respectively, $\Ru$ and  $\Ru'$  in $\End(U\tp W)$.
The same convention is adopted for the matrices $\breve{\Ru}$, $\breve{\Ru'}$ and $\breve{\Sc}$.
By $\Sc_{U}$ we understand the image of the matrix $\Sc^{\chi}$;
then
$\Sc_{U,W}= \breve{\Ru'}_{U,W}\Sc_{U} \breve{\Ru}_{U,W}$.

\section{Representations of $\Ha$ associated with skew Young diagrams}
\label{sYd}
With any skew Young diagram one can associate a one-parameter family of representations
of the Yangian $\Yg(\g\l_N)$. Such representations are called elementary in \cite{NT}.
In fact, elementary representations are obtained from representations of the algebra $\Ha$
via the homomorphism $\Yg(\g\l_N)\to \Ha$. Here we recall the construction.

Consider a non-negative weakly descending sequence $\la=(\la_1,\la_2,\ldots,)$ assuming only a finite
number of non-zero terms. Such a sequence is identified with a Young diagram, which is
a subset in the integer lattice $\Z^2$ specified by the conditions $\{(i,j)\in \Z^2 \>|\>0<i,\> 0< j\leqslant \la_i\}$.
A generally adopted convention is to choose the coordinate system in $\Z^2$ as shown in Figure 1
for the case $\la=(6,5,3,1,0,\ldots)$. The diagram is presented by the set of boxes centered in the lattice nodes.
$$
\begin{picture}(160,120)
\put(-10,110){\vector(0,-1){20}}
\put(-10,110){\vector(1,0){20}}

\put(0,100){\line(0,-1){80}}
\put(0,100){\line(1,0){120}}
\put(0,80){\line(1,0){120}}
\put(0,60){\line(1,0){100}}
\put(0,40){\line(1,0){60}}
\put(0,20){\line(1,0){20}}
\put(20,100){\line(0,-1){80}}
\put(40,100){\line(0,-1){60}}
\put(60,100){\line(0,-1){60}}
\put(80,100){\line(0,-1){40}}
\put(100,100){\line(0,-1){40}}
\put(120,100){\line(0,-1){20}}
\put(-10,80){$i$}
\put(15,106){$j$}
\put(45,){Figure 1.}
\end{picture}
\qquad
\begin{picture}(160,100)
\put(0,60){\line(0,-1){40}}
\put(60,100){\line(1,0){60}}
\put(40,80){\line(1,0){80}}
\put(0,60){\line(1,0){100}}
\put(0,40){\line(1,0){60}}
\put(0,20){\line(1,0){20}}
\put(20,20){\line(0,1){40}}
\put(40,40){\line(0,1){40}}
\put(60,100){\line(0,-1){60}}
\put(80,100){\line(0,-1){40}}
\put(100,100){\line(0,-1){40}}
\put(120,100){\line(0,-1){20}}
\put(50,){Figure 2.}
\end{picture}
$$
We restrict ourself only with diagrams whose columns do not exceed the dimension $N$ of the basic vector space.

Consider  a pair $(\la,\mu)$ of Young diagrams satisfying the condition
$\la_i\geqslant \mu_i$. The skew diagram $\la/\mu$ is the complement
to the diagram $\mu$ in the diagram $\la$.
For example, for  $\la=(6,5,3,1,0,\ldots)$ and $\mu=(3,2,0,\ldots)$
the skew diagram $\la/\mu$ is depicted in Figure 2.

For every skew diagram $\omega=\la/\mu$ we put $|\omega|=\sum_{i}\la_i-\mu_i\in \Z_+$,
the number of boxes in $\omega$.
There are $|\omega| !$ possibilities of filling $\omega$ with the integers $1,\ldots,|\omega|$.
Every such a distribution is called skew Young tableau.
The diagram $\omega$ is called the shape of a tableau.
We shell be dealing with the only tableau of a given shape. It is obtained by filling the boxes
from top to bottom  and from left to right starting from the leftmost column to the right (the so called column standard tableau).
For the skew diagram from Figure 2 the column tableau is displayed in Figure 3.
Thus by a skew diagram we understand the distinguished  skew tableau and use the same notation for it.

The content of a skew tableau is a $\Z$-valued function on the interval $1,\ldots, |\omega|$
defined by $p\mapsto c_p:=j_p-i_p$, where $j_p$ and $i_p$ are respectively the
numbers of column and row where the integer $p$ is allocated.
The tableau of contents for the tableau from Figure 3 is displayed in Figure 4.
$$
\begin{picture}(160,100)
\put(0,60){\line(0,-1){40}}
\put(60,100){\line(1,0){60}}
\put(40,80){\line(1,0){80}}
\put(0,60){\line(1,0){100}}
\put(0,40){\line(1,0){60}}
\put(0,20){\line(1,0){20}}
\put(20,20){\line(0,1){40}}
\put(40,40){\line(0,1){40}}
\put(60,100){\line(0,-1){60}}
\put(80,100){\line(0,-1){40}}
\put(100,100){\line(0,-1){40}}
\put(120,100){\line(0,-1){20}}
\put(8,46){$1$}
\put(8,26){$2$}
\put(28,46){$3$}
\put(48,66){$4$}
\put(48,46){$5$}
\put(68,86){$6$}
\put(68,66){$7$}
\put(88,86){$8$}
\put(88,66){$9$}
\put(104,86){$10$}
\put(50,){Figure 3.}
\end{picture}
\begin{picture}(160,100)
\put(0,60){\line(0,-1){40}}
\put(60,100){\line(1,0){60}}
\put(40,80){\line(1,0){80}}
\put(0,60){\line(1,0){100}}
\put(0,40){\line(1,0){60}}
\put(0,20){\line(1,0){20}}
\put(20,20){\line(0,1){40}}
\put(40,40){\line(0,1){40}}
\put(60,100){\line(0,-1){60}}
\put(80,100){\line(0,-1){40}}
\put(100,100){\line(0,-1){40}}
\put(120,100){\line(0,-1){20}}
\put(3,46){$-2$}
\put(3,26){$-3$}
\put(23,46){$-1$}
\put(48,66){$1$}
\put(48,46){$0$}
\put(68,86){$3$}
\put(68,66){$2$}
\put(88,86){$4$}
\put(88,66){$3$}
\put(108,86){$5$}
\put(50,){Figure 4.}
\end{picture}
$$
With a skew Young diagram we associate a line in the affine space $\C^{|\omega|}$:
$$\C_{\omega}:=\{(v_1,\ldots,v_{|\omega|})\in \C^{|\omega|}\>|\>v_p=z+c_p\},$$
where $z$ is a complex parameter.

Let $\C(z)$ denote the field of rational functions.
Denote by $\eta_\omega$ the ring homomorphism $\C[v_p,v^{-1}_p]\to \C(z)$ corresponding
to the embedding $\C_\omega\hookrightarrow \C^{|\omega|}$.
We shall use the same notation for the  homomorphism $\End(\C^{|\omega|})[v_p,v^{-1}_p]\to \End(\C^{|\omega|})(z)$.

Consider the hyperplane $\Dc_\omega$ in the affine space $\C^{|\omega|}$ defined by the following condition:
$(v_1,\ldots, v_{|\omega|})\in \Dc_\omega$ if and only if $v_p-c_p=v_q-c_q$ once $p$ and $q$ appear in the same column.
Obviously $\C_\omega\subset \Dc_\omega$. It is known,  \cite{C}, that the restriction of the matrix $\breve{R}^{(|\omega|)}:=\bigl(\prod_{p<q}\frac{1}{v_p-v_q}\bigr)R^{(|\omega|)}$
to $\Dc_\omega$ is regular at $(c_1,\ldots,c_{|\omega|})$.  Let $F_\omega \in \End^{\tp |\omega|}(\C^{N})$ denote its limit
as $v_p\to c_p$. As follows from the shift invariance of the Young matrix $R(u,v)$, this limit can be obtained
as $v_p\to c_p+z$ for any $z$. In other words, as $(v_1,\ldots, v_{|\omega|})$ tends to any point on the line $\C_\omega$.

The operator $F_\omega$ is degenerate; in particular, if  $\omega=\la/\mu$ is an ordinary Young diagram
($\mu=\emptyset$), then $F_\omega$ is proportional to the Young idempotent associated with
the corresponding column tableau.

Put $V_\omega$ to be the image of the operator $F_\omega\in\End^{\tp |\omega|}(\C^N)$. Applying $F_\omega$
to the module
$V(v_1)\tilde \tp \ldots \tilde \tp V(v_{|\omega|})$,
where the point $(v_1,\ldots, v_{|\omega|})$ is restricted to the line $\C_\omega$,
we get an $\Ha$-module, denoted further by $V_\omega(z)$.
As a vector space $V_\omega(z)=V_\omega\tp \C_\omega(z)$, where $\C_\omega(z)$ is the image of the homomorphism $\eta_\omega$ in $\C(z)$.

We regard $V_\omega(z)$ as a Yangian module under the homomorphism $\Yg(\g\l_N)\to \Ha$.
When $\omega$ is an ordinary Young diagram, $\V_\omega(z)$ is the so called shifted evaluation $\Yg(\g\l_N)$-module. The
corresponding representation  factors
through the evaluation homomorphism $\Yg(\g\l_N)\to \Ug(\g\l_N)$.
The specialization of the representation $\Yg(\g\l_N)\to \End\bigr(V_\omega(z)\bigr)$ at the point $z\in \C_\omega$ can be factorized in the composition
$$\Yg(\g\l_N)\to \Yg(\g\l_N)\to \End(V_\omega).$$
Here the first arrow is the shift automorphism $T(u)\mapsto T(u-z)$ and the second arrow is the specialization of the representation to
$z=0$.
Note that the latter is well defined for $\Yg(\g\l_N)$.

The subject of our further study is the collection of modules $V_{\omega_1}(z_1)\tilde \tp \ldots \tilde \tp V_{\omega_\ell}(z_\ell)$.
We call them {\em fusion} modules.
Remark that the opposite tensor product $\tilde \tp$ of $\Ha$-modules becomes the tensor product of $\Yg(\g\l_N)$-modules,
as the homomorphism $\Yg(\g\l_N)\to \Ha$ is anti-coalgebra.


\section{Contragredient  diagrams}
\label{Cd}
In the present section we introduce  "conjugate"  skew Young tableaux via a certain
involutive operation on  tableaux. Conjugate tableaux are associated
with contregredient representations of the bialgebra $\Ha$ to be defined later on.

Let $t^{(n)}$ denote the involutive anti-automorphism of the algebra $\End^{\tp n}(\C^N)$ naturally extending the transposition $t$,
that is $t^{(n)}=\tp_{i=1}^n t$. We use the same notation for its extension to the algebra $\Mg^{\tp n}$.
\begin{lemma}
\label{F_t_invariant}
The operator $F_\omega$ is invariant under the involution  $t^{(n)}$, where $n=|\omega|$.
\end{lemma}
\begin{proof}
It suffices to show that the operator $R^{(n)}$ is $t^{(n)}$-invariant.
It is a standard fact from the theory of Hopf algebra twist that the intertwiner
$\Ru^{(n)}$ can be presented
as $(\Delta^{(k)}_{op}\tp \Delta^{(m)}_{op})(\Ru)(\Ru^{(k)}\tp \Ru^{(m)})$, where
$k+m=n$.
This implies that $\Ru^{(n)}$ can be written in either form
$\overrightarrow{\prod}_{p<q}\Ru_{pq}=\Ru^{(n)}=\overleftarrow{\prod}_{p<q}\Ru_{pq}$
(for $n=3$ this is simply the Yang-Baxter equation).
Projecting these equalities to $\Mg^{\tp n}$ and taking into account
$R^t(u,v)=R(u,v)$ we readily prove the assertion.
\end{proof}

Let $\si_n$ be the element of the symmetric group $\Sg_n$ acting by $\si_n(i)=n+1-i$, for $i=1,\ldots, n$.
Denote by $\hat \si_n$ its image in $\End^{\tp }(\C^N)$. Introduce the following involutive automorphisms of the algebra
$\Mg^{\tp n}$ extending the involutions
\be
s^{(n)} (f)(u_i)&:=&f(-u_i),\nn\\
\al^{(n)} (f)(u_i)&:=&f(u_{\si_n(i)}),\nn\\
\bt^{(n)}(A)&:=& \hat\si_n A \hat\si_n\nn
\ee
for $f\in \C[v_i,v_i^{-1}]$ and $A\in \End^{\tp n}(\C^N)$.
The operations $s^{(n)}$, $\al^{(n)}$, $\bt^{(n)}$, and $t^{(n)}$ commute with each other.
Observe that the restriction to $\Mg^{\tp n}$ of the involution $\tau$ coincides with $t^{(n)}s^{(n)}$.

Consider the transposition $(i,j)\mapsto (j,i)$ of the lattice $\Z_2$. Under this transformation,
a Young diagram $\la$ goes over to a diagram $\la^*$. Let $\omega=\la/\mu$ be a skew Young diagram.
Consider the transformation $(i,j)\mapsto (\la_1+1-j,\la_1^*+1-i)$ of the lattice $\Z_2$ (rotation by $180^\circ$ around the center of $\omega$).
Under this transformation, $\omega$ goes over to another skew Young diagram, which we denote by $\omega^\sharp$.
Obviously $(\omega^\sharp)^\sharp=\omega$.
Figure 5 gives such an example.
$$
\begin{picture}(120,100)
\put(-25,60){$\omega=$}
\put(0,60){\line(0,-1){40}}
\put(60,100){\line(1,0){60}}
\put(40,80){\line(1,0){80}}
\put(0,60){\line(1,0){100}}
\put(0,40){\line(1,0){60}}
\put(0,20){\line(1,0){20}}
\put(20,20){\line(0,1){40}}
\put(40,40){\line(0,1){40}}
\put(60,100){\line(0,-1){60}}
\put(80,100){\line(0,-1){40}}
\put(100,100){\line(0,-1){40}}
\put(120,100){\line(0,-1){20}}
\end{picture}
\mbox{Figure 5.}
\begin{picture}(130,100)
\put(-15,60){$\omega^\sharp=$}
\put(0,20){\line(1,0){60}}
\put(0,40){\line(1,0){80}}
\put(20,60){\line(1,0){100}}
\put(60,80){\line(1,0){60}}
\put(100,100){\line(1,0){20}}
\put(0,20){\line(0,1){20}}
\put(20,20){\line(0,1){40}}
\put(40,20){\line(0,1){40}}
\put(60,20){\line(0,1){60}}
\put(80,40){\line(0,1){40}}
\put(100,60){\line(0,1){40}}
\put(120,60){\line(0,1){40}}
\end{picture}
$$
Combined with the inversion $\si_{|\omega|}(i)={|\omega|+1-p}$, $i=1,\ldots, |\omega|$, the mapping $\omega\mapsto \omega^\sharp$
extends to a mapping of tableaux of the corresponding shapes; obviously this preserves the column tableaux.

The  involution  $v_p\mapsto-v_{|\omega|+1-p}$ of the affine space $\C^{|\omega|}$
induces the commutative diagram
\be
\begin{array}{ccccccccc}
\C&\simeq &\C_\omega&\longrightarrow &\Dc_\omega&\longrightarrow &\C^{|\omega|}\\
\downarrow &&\downarrow &&\downarrow &&\downarrow &\\
\C&\simeq &\C_{\omega^\sharp}&\longrightarrow &\Dc_{\omega^\sharp}&\longrightarrow &\C^{|\omega|}\\
\end{array}
\label{sharp}
\ee
It is easy to check that the leftmost downward arrow acts by
\be
\label{shift_inversion}
z\mapsto z^\sharp:=-z-c.
\ee
The integer $c$ is expressed through the contents $c^\omega_p$ and $c^{\omega^\sharp}_p$ of the column tableaux
$\omega$ and $\omega^\sharp$ by
$c=c^\omega_1+c^{\omega^\sharp}_{|\omega|}=c^\omega_p+c^{\omega^\sharp}_{\si(p)}$. Here
$\si=\si_{|\omega|}$; thus the right-hand side is independent on $p$.

Let $\eta'_{\omega}$ denote the homomorphism from the ring $\C[v_i,v_i^{-1}]$ to the field of rational functions
on $\Dc_\omega$ corresponding to the embedding $\Dc_\omega \hookrightarrow \C^{|\omega|}$.
The rightmost square of the diagram (\ref{sharp}) can be written
as $\theta \eta'_{\omega^\sharp}=\eta'_{\omega} \al^{(n)}s^{(n)}$, where
$\theta \colon \C(\Dc_{\omega^\sharp})\to \C(\Dc_{\omega})$ is an
isomorphism of the corresponding functional rings.
Our next goal is to relate the intertwiners $F_{\omega^\sharp}$ and $F_{\omega}$.
\begin{propn}
Let $\omega$ be a skew Young diagram with $n=|\omega|$ boxes.
Then $F_{\omega^\sharp}=\hat\si_n F_{\omega}\hat\si_n$.
\end{propn}

\begin{proof}
As we argued in the proof of Lemma \ref{F_t_invariant},
 we can reverse the order of the $\Ru_{pq}$-factors in the definition of the operator $\Ru^{(n)}$.
Thus we have for $(\theta\eta'_{\omega^\sharp})(R^{(n)})$:
\be
&&(\theta\eta'_{\omega^\sharp})\bigl(\overleftarrow{\prod_{p<q}}R_{p,q}(v_p,v_q)\bigr)
=\eta'_\omega\al^{(n)}s^{(n)}\bigl(\overleftarrow{\prod_{p<q}}R_{p,q}(v_p,v_q)\bigr)
=\eta'_\omega\al^{(n)}\bigl(\overleftarrow{\prod_{p<q}}R_{p,q}(v_q,v_p)\bigr)
\nn\\
&&=
\eta'_\omega\bigl(\overleftarrow{\prod}_{p<q}R_{p,q}(v_{\si_n(q)},v_{\si_n(p)})\bigr)=\eta'_\omega\bigl(\overleftarrow{\prod}_{p>q}R_{\si_n(p),\si_n(q)}(v_q,v_p)\bigr),
\nn
\ee
where we have made the substitution $p\mapsto \si_n(p),q\mapsto \si_n(q)$.
Using  the equalities $R_{ij}(u,v)=R_{ji}(u,v)$ and $\overleftarrow{\prod}_{p>q}=\overrightarrow{\prod}_{q<p}$,
we rewrite the last expression in the form
$$
\eta'_\omega\bigl(\overrightarrow{\prod}_{q<p}R_{\si_n(p),\si_n(q)}(v_q,v_p)\bigr)=\hat\si_n\eta'_\omega\bigl(R^{(n)}\bigr)\hat\si_n.
$$
Dividing the resulting equality by $\prod_{p<q}(v_p-v_q)$ we get
\be
\theta\bigl(\breve{R}^{(n)}|_{\Dc_{\omega^\sharp}}\bigr)&=&\hat\si_n\breve{R}^{(n)}|_{\Dc_\omega}\hat\si_n.
\label{eq:aux1}
\ee
Suppose that $d:=(v_1,\ldots,v_{|\omega|})\in \Dc_\omega\subset \C^{\tp |\omega|}$.
Let $\theta_*$ denote the isomorphism $\Dc_\omega\to \Dc_{\omega^\sharp}$ from the diagram (\ref{sharp}).
As follows from middle square of the diagram (\ref{sharp}), the point $\theta_*(d)\in \Dc_{\omega^\sharp}$ tends to a point on the line $\C_{\omega^\sharp}$
as $d\in \Dc_\omega$ tends to a point on the line $\C_{\omega}$.
Taking the  limit $v_p\to c_p$  of the both sides
of the equation (\ref{eq:aux1}) we prove the proposition.
\end{proof}

\section{Contragredient $\Ha$-modules}
\label{CHm}
Denote by $V_\omega^\sharp(z)$ the $\Ha$-module $V_{\omega^\sharp}(z^\sharp)$ and by $\rho^\sharp_{\omega}$ the corresponding representation  homomorphism
$\Ha\to \End(V_{\omega}^\sharp)(z)$.
The module $V_\omega^\sharp(z)$ can be obtained as a submodule in the tensor product of $n$-copies of the defining representations
upon the specialization $\eta^\sharp_{\omega}=\eta_{\omega}\al^{(n)} s^{(n)}\colon \C[z_p,z^{-1}_p]\to \C_{\omega^\sharp}(z)$.
The homomorphism $\eta^\sharp_{\omega}$ can also be presented as the composition of the specialization $\eta_{\omega^\sharp}$
and the subsequent involution $z\mapsto z^\sharp$ of the line $\C_{\omega^\sharp}$, cf. the diagram (\ref{sharp}).

For any $\Ha$-module $(W,\rho)$ define the contragredient representation $\rho^*$
on the dual vector space $W^*$ by setting $\rho^*(h)(\phi)(w):=\phi\bigl((\rho\circ\tau)(h)(w)\bigr)$
for all $h\in \Ha$, $\phi\in W^*$, and $w\in W$. Equivalently, put $\rho^*=t\rho\tau$.
It is easy to see that $(U\tilde\tp W)^*=U^*\tilde\tp W^*$ for any pair of modules $U$ and $W$ (recall that $\tau$ preserves
the comultiplication).

\begin{propn}
\label{duality}
Let $\omega$ be a skew Young diagram.
Then the module $V_\omega^\sharp(z)$ is isomorphic to $V_\omega^*(z)$. Namely, for all $h\in \Ha$:
$$
\rho_{\omega}^\sharp(h)F_{\omega^\sharp}=\hat\si_n\bigl((\rho_{\omega}\circ\tau )(h)F_{\omega}\bigr)^{t^{(n)}}\hat\si_n=\hat\si_nF_{\omega}\bigl((\rho_{\omega}\circ\tau) (h)\bigr)^{t^{(n)}}\hat\si_n.
$$
\end{propn}

\begin{proof}
The involution $\al^{(n)}\bt^{(n)}$ is the total permutation of $\Mg^{\tp n}$ reversing the order of the tensor factors.
It is sufficient to restrict the consideration to  $h\in \Mg^{\tp n}\cap \Ha$. Then
for all $\omega$ we have
$$
\rho_{\omega}(h)F_{\omega}= \bigl(\eta_\omega\al^{(n)}\bt^{(n)}\bigr)(h) F_\omega=F_\omega\eta_\omega( h) ,
$$
which is a corollary of (\ref{intertwining}).
Hence
$$
(\rho_{\omega}^\sharp\tau)(h)F_{\omega^\sharp}
=F_{\omega^\sharp}\bigl(\eta_{\omega}^\sharp \tau\bigr) (h)
=F_{\omega^\sharp}\bigl(\eta_{\omega}^\sharp t^{(n)}s^{(n)}\bigr) (h),
$$
where $\eta_\omega^\sharp$ is the composition of $\eta_{\omega^\sharp}$ and the involution $z\mapsto z^\sharp$.
It can also be presented as $\eta_{\omega}^\sharp=\eta_{\omega}\al^{(n)} s^{(n)}$, as follows from the diagram (\ref{sharp}).
We continue the calculation equality as
$$
=F_{\omega^\sharp}\bigl(\eta_{\omega}\al^{(n)} s^{(n)}t^{(n)}s^{(n)}\bigr) (h)=F_{\omega^\sharp}\bigl(\eta_{\omega}\al^{(n)}t^{(n)}\bigr) (h)
=
t^{(n)}\Bigl(\bigl(\eta_{\omega}\al^{(n)}\bigr) (h)t^{(n)}(F_{\omega^\sharp})\Bigr).
$$
Now  use Lemma \ref{F_t_invariant} and  substitute $F_{\omega^\sharp}=\hat\si_n F_{\omega}\hat\si_n=(t^{(n)}\bt^{(n)})(F_\omega)$.
We continue the calculation by
$$
=
t^{(n)}\Bigl(\bigl(\eta_{\omega}\al^{(n)}\bigr) (h)\bt^{(n)}(F_{\omega})\Bigr)
=
t^{(n)}\bt^{(n)}\Bigl(\bigl(\eta_{\omega}\bt^{(n)}\al^{(n)}\bigr) (h)F_{\omega})\Bigr)
=
\hat\si_n\Bigl(\rho_{\omega}(h)F_{\omega}\Bigr)^{t^{(n)}}\hat\si_n.
$$
The proof is complete.
\end{proof}

We extend the operation $\sharp$ over all fusion modules. Specifically,
for any finite set of skew Young diagrams $\{\omega_i\}_{i=1}^\ell$ we put
$$
\Bigr(V_{\omega_1}(z_1)\tilde \tp \ldots \tilde \tp V_{\omega_\ell}(z_\ell)\Bigr)^\sharp
:=
V^\sharp_{\omega_1}(z_1)\tilde \tp \ldots \tilde \tp V^\sharp_{\omega_\ell}(z_\ell)
=
V_{\omega_1^\sharp}(z_1^\sharp)\tilde \tp \ldots \tilde \tp V_{\omega_\ell^\sharp}(z_\ell^\sharp).
$$
Note that the involution $z_i\mapsto z_i^\sharp$ is defined via the skew diagram $\omega_i$ for each $i=1,\ldots, \ell$
by the  formula (\ref{shift_inversion}).
By Proposition \ref{duality}, the module $X^\sharp$ is isomorphic to the contragredient module $X^*$
for any fusion module $X$.
It is easy to see that $(X^\sharp)^\sharp=X$.

For any pair of skew Young diagrams $\vt,\omega$, the matrix $\Ru_{V_{\vt},V_\omega}(z_1,z_2)$
is a rational function of the variables $z_1,z_2$ with values in $\End({V_\vt\tp V_\omega})$.
It is known that the matrix $\Ru_{V_{\vt},V_\omega}(z_1,z_2)$ is invertible if the point
$(z_1,z_2)$ lies outside a finite union of lines in $\C^2$. We need an analogous statement for
the matrices $\Ru_{V_{\omega}^\sharp,V_\omega}(z,z)=\Ru_{V_{\omega^\sharp},V_\omega}(z^\sharp,z)$
and   $\Sc_{V_\omega}(z)$ (for the definition of $\Sc_{V_\omega}$ see Section \ref{Sec:H-mod}).
\begin{lemma}
\label{invertibility}
The matrices $\Ru_{V_{\omega^\sharp},V_\omega}(z^\sharp,z)$ and  $\Sc_{V_\omega}(z)$ are invertible except for a finite number of points of $\C$.
\end{lemma}
\begin{proof}
The matrix $\Ru_{V_{\omega^\sharp},V_\omega}(z^\sharp,z)$ can be presented explicitly as
$$
\Ru_{V_{\omega^\sharp},V_\omega}(z^\sharp,z)=\prod_{i=n}^1\prod_{i=1}^n \Bigr(1+\frac{P_{i j}}{2z+c_{\si(i)}+c_j}\Bigl)|_{V_{\omega^\sharp}\tp V_{\omega}},
$$
where $\si=\si_{|\omega|}\in \Sg_{|\omega|}$ is the permutation reversing the order of sequence $1,\ldots,|\omega|$.
Here the indices $i$ and $j$ refer to the first and the second tensor factors, respectively.
The order of factors in the ordered products here and further on is determined by the initial and final values of the indices
are assumed to be directed from left to right.

If follows from the above presentation that the singular points of the operator valued function $\Ru_{V_\omega^\sharp,V_\omega}(z^\sharp,z)$ are contained in
$\{-\frac{c_{\si(i)}+c_j}{2},-\frac{c_{\si(i)}+c_j\pm 1}{2}\}_{i<j}$.
This proves the statement for $\Ru_{V_\omega^\sharp,V_\omega}(z^\sharp,z)$.

The equation (\ref{char_eq}) yields the following expression for the matrix $\Sc_{V_\omega}(z)$:
$$
\Sc_{V_\omega}(z)=(-1)^{\frac{|\omega|(|\omega|-1)}{2}}\overleftarrow{\prod_{p<q}}\bigr(2z+c_p+c_q+Q_{p,q}\bigl)|_{V_{\omega}}.
$$
The operator $Q\in \End^{\tp 2}(\C^N)$ is defined in Section \ref{YtY} and obeys the equality $Q^2=NQ$.
The matrix $\Sc_{V_\omega}(z)$  is certainly invertible outside the finite set $\{-\frac{c_p+c_q}{2},-\frac{c_p+c_q+N}{2}\}_{p<q}$.
This completes the proof.
\end{proof}
Recall once again that we use the opposite coproduct in $\Ha$ to  construction the tensor product of modules;
that accounts for the reversed order of the factors in the formulas for and $\Ru_{V_{\omega^\sharp},V_\omega}$ and $\Sc_{V_\omega}$.
\section{Some explicit factorization formulas}
\label{Sec:factorization}
In this section we write down explicit formulas for some operators that we use in what follows.
Put $W=V_{\vt_1}(w_1)\tilde\tp \ldots\tilde\tp V_{\vt_k}(w_k)$ and $Z=V_{\omega_1}(z_1)\tilde\tp \ldots\tilde\tp V_{\omega_\ell}(z_\ell)$,
where $\{\vt_i\}_{i=1}^k$ and $\{\omega_j\}_{j=1}^\ell$ are two families of skew Young diagrams.
Put $m_i:=|\vt_i|$ and $n_j:=|\omega_j|$ for $i=1,\ldots, k$ and $j=1\ldots, \ell$;
denote  also $m:=\sum_{i=1}^k m_i$ and $n:=\sum_{j=1}^\ell m_j$.

As vector spaces, $W$ and $V$ are the images of the operators $\tp_{i=1}^k F_{\vt_i}\in \End^{\tp m}(\C^N)$ and $\tp_{j=1}^\ell F_{\omega_j}\in \End^{\tp m}(\C^N)$, respectively.
By definition,
$$
(\C^N)^{\tp m}\tp (\C^N)^{\tp n}=\Bigl((\C^N)^{\tp m_1}\tp\ldots\tp  (\C^N)^{\tp m_k}\Bigr)\tp \Bigl((\C^N)^{\tp n_1}\tp\ldots\tp  (\C^N)^{\tp n_\ell}\Bigr).
$$
Using the identities $(\Delta_{op}\tp \id)(\Ru)=\Ru_{23}\Ru_{13}$, $(\id\tp \Delta_{op})(\Ru)=\Ru_{12}\Ru_{13}$
with the opposite comultiplication  we can write
\be
\Ru_{W,Z}(w_1,\ldots, w_k, z_1,\ldots, z_\ell)
&=&\prod_{i=k}^1\prod_{j=1}^\ell\Ru_{V_{\vt_i},V_{\omega_j}}(w_i,z_j)\nn,\\
\Ru'_{W,Z}(w_1,\ldots, w_k, z_1,\ldots, z_\ell)
&=&\prod_{i=k}^1\prod_{j=\ell}^1\Ru'_{V_{\vt_i},V_{\omega_j}}(w_i,z_j)\nn.
\nn
\ee
Here the index $i=1,\ldots, k$ enumerates the factors $(\C^N)^{\tp m_i}$ in  $(\C^N)^{\tp m}$,
while the index $j=1,\ldots, \ell$ enumerates the factors $(\C^N)^{\tp n_j}$ in  $(\C^N)^{\tp n}$.
Explicitly we have
\be
\Ru_{V_{\vt_i},V_{\omega_j}}(w_i,z_j)
&=&\prod_{p_i=m_k}^1\prod_{q_j=1}^{n_\ell} (u_{p_i}-v_{q_j}-P_{p_i,q_j})|_{V_{\vt_i}\tp V_{\omega_j}},\nn\\
\Ru'_{V_{\vt_i},V_{\omega_j}}(w_i,z_j)
&=&\prod_{p_i=m_k}^1\prod_{q_j=n_\ell}^{1} (-u_{p_i}-v_{q_j}-Q_{p_i,q_j})|_{V_{\vt_i}\tp V_{\omega_j}}.\nn
\nn
\ee
Here
$
u_{p_i}=w_i+c_{p_i}$ and $v_{q_j}=z_j+c_{q_j}$,
where $c_{p_i}$ denotes the contents of the tableau $\vt_i$ while $c_{q_j}$ the contents of the tableau $\omega_j$.
The index $p_i=1,\ldots, m_i$
refers to the $\C^N$-factor in $(\C^N)^{\tp m_i}$ and, similarly, the index $q_j=1,\ldots, n_j$ refers to the $\C^N$-factor in $(\C^N)^{\tp n_j}$.

Specializing to the case $W=Z^\sharp$, we obtain
\be
\Ru_{Z^\sharp,Z}(z_1,\ldots, z_\ell,z_1,\ldots, z_\ell)
&=&\prod_{i=k}^1\prod_{j=1}^\ell\Ru_{V_{\omega_i}^\sharp,V_{\omega_j}}(z_i,z_j)\nn,\\
\Ru'_{Z^\sharp,Z}(z_1,\ldots, z_\ell,z_1,\ldots, z_\ell)
&=&\prod_{i=k}^1\prod_{j=\ell}^1\Ru'_{V_{\omega_i}^\sharp,V_{\omega_j}}(z_i,z_j)\nn,
\nn
\ee
where
\be
\Ru_{V_{\omega_i}^\sharp,V_{\omega_j}}(z_i,z_j)
&=&\Ru_{V_{\omega^\sharp_i},V_{\omega_j}}(z^\sharp_i,z_j)=
\prod_{p_i=n_\ell}^1\prod_{q_j=1}^{n_\ell} (-v_{\si(p_i)}-v_{q_j}-P_{p_i,q_j}) |_{V_{\omega_i^\sharp}\tp V_{\omega_j}},\nn\\
\Ru'_{V_{\omega^\sharp_i},V_{\omega_j}}(z_i,z_j)
&=&
\Ru_{V_{\omega_i}^\sharp,V_{\omega_j}}(z^\sharp_i,z_j)'=
\prod_{p_i=n_\ell}^1\prod_{q_j=n_\ell}^{1} (u_{\si(p_i)}-v_{q_j}-Q_{p_i,q_j})|_{V_{\omega_i^\sharp}\tp V_{\omega_j}}.\nn
\nn
\ee
Equation (\ref{char_eq}) descends to the formula
\be
\Sc_Z(z_1,\ldots, z_\ell)&=&\prod_{i=\ell}^1\Bigl(\Sc_{V_{\omega_i}}(z_i)\prod_{j=i-1}^1\Ru'_{V_{\omega_i},V_{\omega_j}}(z_i,z_j)\Bigr),
\nn
\quad \mbox{where}\\
\Sc_{V_{\omega_i}}(z_i)&=&\prod_{p_i=n_i}^1\Bigl(\prod_{q_j=p_i-1}^1R'_{p_i,q_j}(v_{p_i},v_{q_j})\Bigr)|_{V_{\omega_i}},
\nn
\ee
see the proof of Lemma \ref{invertibility}.
As earlier, we assume $v_{q_i}=z_i+c_{q_i}$,
where  $c_{q_i}$ denote the contents of the tableau $\omega_i$.
Recall that we use the opposite coproduct of $\Ha$ in definition of the tensor product of modules.

The matrices $\breve{\Ru}_{W,Z}(w_1,\ldots, w_k, z_1,\ldots, z_\ell)$  and $\breve{\Ru}_{W,Z}'(w_1,\ldots, w_k, z_1,\ldots, z_\ell)$ are expressed
by the same formulas, where $R_{p_i,q_j}(u_{p_i},v_{q_j})$ and $R'_{p_i,q_j}(u_{p_i},v_{q_j})$ are replaced by $\breve{R}_{p_i,q_j}(u_{p_i},v_{q_j})$
and $\breve{R}'_{p_i,q_j}(u_{p_i},v_{q_j})$, respectively.

\section{Irreducibility of fusion modules over $\Yg(\g_N)$ at generic point}
\label{Irreducibility_gen}
In this section we prove that fusion modules are irreducible over the algebra  $\Xg(\g_N)$
at generic point. That will also imply irreducibility over the twisted Yangian $\Yg(\g_N)$, as the
representations in question factor through the projection $\Xg(\g_N)\to\Yg(\g_N)$.
We rely on an observation providing a
necessary condition for modules over associative algebras to be irreducible.
The method was applied in \cite{NT}
to Yangian modules and  in \cite{MN} to twisted Yangian modules associated with diagrams of one row or one column.

We begin  with formulation of the method for an arbitrary associative algebra, call it $\A$.
Fix a finite dimensional $\A$-module  $(Z,\rho)$  and
a finite dimensional vector space $W$.  Identify elements of
$\End(W)\tp \End(Z)$  with linear maps from $\End(W)$ to $ \End(Z)$
using pairing via the trace $\Tr_W$.
Let  $\dot{U}$ be a punctured neighborhood of $0$ in $\C$. Suppose
$\Phi\colon\dot{U}\to \End(W)\tp \rho(\A)$ is an analytical function
with the Laurent decomposition $\Phi(\zeta)=\zeta^{r}\Phi_0+O(\zeta^{r+1})$, $r\in \Z$.
\begin{propn}
\label{Prop_irred}
Suppose that the leading  Laurent coefficient $\Phi_0$ of the operator function $\Phi(\zeta)$ at $\zeta=0$
implements a surjective mapping $\End(W)\to \End(Z)$.
Then the $\A$-module $Z$ is irreducible.
\end{propn}
\begin{proof}
For any $A\in \End(W)$ and all $\zeta$ from $\dot{U}$ the matrix
$\zeta^{-r}\Phi(\zeta)(A)$ belongs to $\rho(\A)$.
Therefore its limit $\Phi_0(A)$ as $\zeta\to 0$ belongs to $\rho(\A)$ too, since the latter
is closed in $\End(Z)$ in the euclidian topology (as a finite dimensional vector space).
Thus $\im \Phi_0=\End(Z)\subset \rho(\A)$ and hence $\End(Z)= \rho(\A)$. This implies irreducibility of $Z$.
\end{proof}
We call $W$ the auxiliary vector space,  $\Phi(\zeta)$ the auxiliary operator function, and its leading term $\Phi_0$
simply the auxiliary operator.

We apply Proposition \ref{Prop_irred} to  the $\Xg(\g_N)$-module $Z=V_{\omega_1}(z_1)\tilde\tp \ldots\tilde\tp V_{\omega_\ell}(z_\ell)$
specialized
at generic point  $(z_1,\ldots, z_\ell)$. The underlying vector space of the module  $Z$ is
the tensor product $V_{\omega_1}\tp \ldots\tp V_{\omega_\ell}$.
We take for $W$ the same vector space $V_{\omega_1}\tp \ldots\tp V_{\omega_\ell}$
regarded as an $\Ha$-module
$V_{\omega_1}(z_1+\zt)\tilde\tp \ldots\tilde\tp V_{\omega_\ell}(z_\ell+\zt)$.
It follows that  $W=Z$ at the point $\zt=0$.

For a fixed point $(z_1,\ldots, z_\ell)$ regard the matrix $\Sc_{W,Z}$ as a rational function of the variable $\zt$ with
values in endomorphisms of the vector space $(V_{\omega_1}\tp \ldots \tp V_{\omega_\ell})^{\tp 2}$.
The matrix $\Sc_{W,Z}$ lies in $\End(W)\tp \rho_Z\bigl(\Xg(\g_N)\bigr)$,
where
$\rho_Z\colon \Xg(\g_N)\to \End(Z)$ is the representation.
We take $\Sc_{W,Z}$ as the auxiliary operator function $\Phi(\zt)$.

\begin{thm}
The $\Xg(\g_N)$-module $V_{\omega_1}(z_1)\tilde\tp \ldots\tilde\tp V_{\omega_\ell}(z_\ell)$ is irreducible
for generic $(z_1,\ldots, z_\ell)$.
\end{thm}
\begin{proof}
Explicitly, the matrix $\Sc_{W,Z}$ is presented as
\be
\label{factorization_S}
\Sc_{W,Z}=\breve{\Ru}'_{W,Z}\Sc_W\breve{\Ru}_{W,Z},
\ee
see Section \ref{Sec:H-mod}.
Let $n_i$ be the number of boxes in the skew diagram $\omega_i$,
so that every $V_{\omega_i}$ is a vector subspace
in $(\C^N)^{\tp n_i}$. Then  $V_{\omega_1}\tp \ldots\tp V_{\omega_\ell}$ is a subspace in $(\C^N)^{\tp n}$, with $n:=\sum_{i=1}^\ell n_i$.

Denote by $P_{(n),(n)}$ the flip of tensor factors in $(\C^N)^{\tp n}\tp (\C^N)^{\tp n}$.
As a rational function of $\zt$, the matrix $\breve{\Ru}_{W,Z}$ has a pole at $\zt=0$. The leading term of the Laurent expansion is equal to
$h_ZP_{(n),(n)}$, where  $h_Z$ is a rational function of $z_1,\ldots, z_\ell$, see \cite{NT}.
Observe that the operator $A\mapsto (\Tr\tp \id)\Bigl((A\tp 1)P_{(n),(n)}\Bigr)$ is identical on $\End^{\tp n}(\C^N)$.

Consider the assignment $\Phi(\zt)\colon A\mapsto (\Tr\tp \id)\bigl((A\tp 1)(\Sc_{W,Z})\bigr)$,
$A\in \End(V_{\omega_1}\tp \ldots\tp V_{\omega_\ell})$,
as an auxiliary operator.
The factorization formula (\ref{factorization_S})
yields the leading Laurent term of the operator $\Phi(\zt)$ at the point $\zt=0$:
\be
\label{auxiliary}
\Phi_0\colon A\mapsto h_Z(\breve{\Ru}'_2)_Z A (\breve{\Ru}'_1)_Z \Sc_Z,
\ee
provided $h_Z$, $\Ru_{Z^\sharp,Z}$, and $\Sc_Z$ are regular and do not vanish at $(z_1,\ldots, z_\ell)$. Here we have used the symbolic notation with
suppressed summation
$\breve{\Ru}'=\breve{\Ru}'_1\tp \breve{\Ru}'_2=\tau(\breve{\Ru}_1)\tp \breve{\Ru}_2$.

By Proposition \ref{duality}, $\breve{\Ru}'_{Z,Z}=(\si'_n\tp 1)(t^{(n)}\tp \id)(\breve{\Ru}_{Z^\sharp,Z})(\si'_n\tp 1)$,
where $\hat \si'_n:=\tp_{i=1}^\ell \hat\si_{n_i}$.
Recall that $Z^\sharp$ is the $\Xg(\g_N)$-module  $V_{\omega^\sharp_1}(z^\sharp_1)\tp \ldots\tp V_{\omega^\sharp_\ell}(z^\sharp_\ell)$.
Hence the map  (\ref{auxiliary}) can be presented as
\be
\Phi_0\colon  A\mapsto h_Z\Psi_2A\>t^{(n)}(\Psi_1)\Sc_Z,
 \nn
\ee
where $\Psi:=(\hat \si'_n\tp 1)(\breve{\Ru}_{Z^\sharp,Z})(\hat \si'_n\tp 1)\in \End^{\tp 2}(V_{\omega_1}\tp \ldots\tp V_{\omega_\ell})$. The operator $\breve{\Ru}_{Z^\sharp,Z}$
is a rational function of  $(z_1,\ldots, z_\ell)$, as well as the operator  $\Sc_Z$.
It follows from Proposition \ref{invertibility} and the factorization formulas of Section \ref{Sec:factorization},  the operators
$h_Z\Sc_Z$ and $\breve{\Ru}_{Z^\sharp, Z}$
are invertible outside an algebraic set in the parameter space $\C^\ell$.
As the transposition $t^{(n)}$ is an anti-algebra automorphism, the inverse to the map  (\ref{auxiliary}) is given
by the assignment
\be
\label{auxiliary_inv}
\Phi_0^{-1}\colon   A\mapsto \frac{1}{h_Z}\bar \Psi_2A\Sc_Z^{-1}t^{(n)}(\bar\Psi_1)
 \nn
 \ee
with $\bar \Psi:=\Psi^{-1}$.
Therefore the mapping
(\ref{auxiliary}) is invertible and hence surjective at generic point.
This proves the theorem.
\end{proof}
Factorization formulas from Section \ref{Sec:factorization} together with Lemma \ref{invertibility}
allow to give a rough estimate
of the set of points where the $\Yg(\g_N)$-module $V_{\omega_1}(z_1)\tilde\tp \ldots\tilde\tp V_{\omega_\ell}(z_\ell)$
is reducible. That is a finite union of hyperplains in $\C^\ell$, each being of the following three forms:
\be
z_i&=&k_i\in \frac{1}{2}\Z,
\label{single}\\
z_i-z_j&=&k^-_{ij}\in \Z,
\label{Yangian}\\
z_i+z_j&=&k^+_{ij}\in \Z.
\label{twYangian}
\ee
The condition (\ref{Yangian}) descends from the function $h_Z$ and it is responsible for reducibility of each pair $V_{\omega_i}(z_i)\tilde\tp V_{\omega_j}(z_j)$
over the Yangian $\Yg(\g\l_N)$, see \cite{NT} for details. The conditions (\ref{single}) reflects a possibility for the $\Yg(\g\l_N)$-module $V_{\omega_i}(z_i)$ to turn reducible upon
restriction to $\Yg(\g_N)$. If all the parameters $z_i$ are not half-integers and $z_i\pm z_j\not \in \Z$ for all pairwise distinct $i,j$, then
the module $V_{\omega_1}(z_1)\tilde\tp \ldots\tilde\tp V_{\omega_\ell}(z_\ell)$ is irreducible over the
twisted Yangian $\Yg(\g_N)$.

\end{document}